\documentclass[10pt]{amsart}

\usepackage{amssymb,amsmath,amsthm}
\usepackage[all]{xy}
\usepackage{eucal}

\theoremstyle{plain}

\newtheorem{prop}[subsection]{Proposition}

\newtheorem{innercustomprop}{Proposition}
\newenvironment{customprop}[1]
  {\renewcommand\theinnercustomprop{#1}\innercustomprop}
  {\endinnercustomprop}
  
\theoremstyle{definition}
\newtheorem{defn}[subsection]{Definition}

\theoremstyle{remark}
\newtheorem{rem}[subsection]{Remark}

\makeatletter 
\let\c@equation\c@subsection 
 
\makeatother

\newcommand{\I}{{ \mathsf{I} }}

\newcommand{\sSet}{{ \mathsf{S} }}

\newcommand{\Spectra}{{ \mathsf{Sp}^\Sigma }}

\newcommand{\SymSeq}{{ \mathsf{SymSeq} }}

\newcommand{\Z}{{ \mathbb{Z} }}

\newcommand{\capA}{{ \mathcal{A} }}
\newcommand{\capB}{{ \mathcal{B} }}
\newcommand{\capX}{{ \mathcal{X} }}

\newcommand{\id}{{ \mathrm{id} }}
\newcommand{\op}{{ \mathrm{op} }}

\newcommand{\Smash}{{ \,\wedge\, }}

\newcommand{\tensor}{{ \otimes }}

\newcommand{\tensorcheck}{{ \check{\tensor} }}

\newcommand{\wequiv}{{ \ \simeq \ }}

\newcommand{\Iso}{{  \ \cong \ }}

\newcommand{\rarrow}{{ \rightarrow }}

\newcommand{\functor}[3]{{ {#1}\colon\thinspace{#2}\rarrow{#3} }}
\newcommand{\function}[3]{{ {#1}\colon\thinspace{#2}\rarrow{#3} }}

\newcommand{\subsetof}{{ \ \subset\ }}

\newdir{ >}{{}*!/-8pt/@{>}}
\newcommand{\monic}{ \ar@{ >->} }

\DeclareMathOperator*{\colim}{colim}

\DeclareMathOperator{\BAR}{Bar}

\DeclareMathOperator{\orbit}{orbit}

\title[Corrigendum to ``Homotopy theory of modules over operads'']{Corrigendum to ``Homotopy theory of modules over operads in symmetric spectra''}

\author{John E. Harper}

\address{Department of Mathematics, The Ohio State University, Newark, 1179 University Dr, Newark, OH 43055, USA}

\email{harper.903@osu.edu}

\begin{document}
\maketitle

\begin{abstract}
Dmitri Pavlov and Jakob Scholbach have pointed out that part of Proposition 6.3, and hence Proposition 4.28(a), of \cite{Harper_Spectra} are incorrect as stated. While all of the main results of that paper remain unchanged, this necessitates modifications to the statements and proofs of a few technical propositions.
\end{abstract}

\section{Introduction}


The author would like to thank Dmitri Pavlov and Jakob Scholbach for pointing out that the description of the cofibrations in the last sentence of Proposition 6.3 of \cite{Harper_Spectra} is incorrect as stated; in general, to verify that a map is a cofibration, it is not enough to be a monomorphism such that $\Sigma_r^\op \times G$ acts freely on the simplices of the codomain not in the image. 

It is well known that the cofibrations in $\sSet_*^G$, equipped with the projective model structure, are precisely the monomorphisms such that $G$ acts freely on the simplices of the codomain not in the image. One way to verify this is to (i) argue that the image of such a map is a subcomplex of the codomain (i.e, the codomain can be built from the image by attaching $G$-cells), and (ii) note that every monomorphism is isomorphic to its image, hence verifying such maps are cofibrations, (iii) conversely, to note that every generating cofibration is such a map, and (iv) hence conclude that every cofibration is such a map, by using the fact that every cofibration is a retract of a (possibly transfinite) composition of pushouts of the generating cofibrations. The problem with our argument for the cofibration description in \cite[6.3]{Harper_Spectra} was a cavalier application of the subcomplex argument (i) above; we ignored the fact that $\Sigma_r^\op\times G$ and $\Sigma_n$ might not act independently. Pavlov and Scholbach kindly pointed out this problem to the author, together with a helpful counterexample to focus one's attention. At the time they were working to generalize the main results in \cite{Harper_Spectra} to  motivic settings (including Hornbostel's results \cite{Hornbostel}, see Remark \ref{rem:Hornbostel}). Their efforts have now appeared in Pavlov and Scholbach \cite{Pavlov_Scholbach}; included in \cite[A]{Pavlov_Scholbach} is their helpful counterexample, together with further discussion related to these cofibrations. 

The following proposition corresponds to the corrected version of \cite[6.3]{Harper_Spectra}.

\begin{customprop}{6.3*}
\label{prop:six_point_three_star}
Let $G$ be a finite group and consider any $n,r\geq 0$. The diagram category $\bigl(\sSet_*^{\Sigma_n}\bigr)^{\Sigma_r^\op\times G}$ inherits a corresponding projective model structure from the mixed $\Sigma_n$-equivariant model structure on $\sSet_*^{\Sigma_n}$. The weak equivalences (resp. fibrations) are the underlying weak equivalences (resp. fibrations) in $\sSet_*^{\Sigma_n}$.
\end{customprop}

The consequence of this misunderstanding of the cofibrations in \cite[6.3]{Harper_Spectra} is that Proposition 4.28(a) of \cite{Harper_Spectra} is incorrect as stated. While all of the main results of that paper remain unchanged, this necessitates modifications to the statements and proofs of a few technical propositions. 

\begin{rem}
\label{rem:Hornbostel}
This corrigendum also applies to the proof of the motivic generalization of our results provided by Hornbostel, namely \cite[3.6, 3.10, 3.15]{Hornbostel}.
\end{rem}

The following proposition corresponds to the corrected version of \cite[4.28]{Harper_Spectra}. For a useful study of additional properties associated to tensor powers of cofibrations, see Pereira \cite{Pereira_spectral_operad}, and more recently, Pavlov and Scholbach \cite{Pavlov_Scholbach}.

\begin{customprop}{4.28*}
\label{prop:four_point_two_eight_star}
Let $B\in\SymSeq^{\Sigma_t^\op}$, $t\geq 1$, and $r,n\geq 0$. If $\function{i}{X}{Y}$ is a cofibration between cofibrant objects in $\SymSeq$ with the positive flat stable model structure, then 
\begin{itemize}
\item [(a)] the map $B\tensorcheck X^{\tensorcheck t}\rarrow B\tensorcheck Y^{\tensorcheck t}$, after evaluation at $[\mathbf{r}]_n$, is a cofibration in $\sSet_*^{\Sigma_t}$ with the projective model structure inherited from $\sSet_*$,
\item [(b)] the map $B\tensorcheck_{\Sigma_t}Q_{t-1}^t\rarrow B\tensorcheck_{\Sigma_t}Y^{\tensorcheck t}$ is a monomorphism.
\end{itemize}
\end{customprop}

Since Propositions 4.29 and 6.11 of \cite{Harper_Spectra} are no longer immediately applicable, we include below the closely related Propositions \ref{prop:four_point_two_nine_star} and \ref{prop:six_point_one_one_star} which describe the technical properties that are actually used in the proofs of the main results in \cite{Harper_Spectra}.

\begin{customprop}{4.29*}
\label{prop:four_point_two_nine_star}
Let $t\geq 1$ and consider $\SymSeq$ and $\SymSeq^{\Sigma_t^\op}$ each with the positive flat stable model structure.
\begin{itemize}
\item [(a)] If $B\in\SymSeq^{\Sigma_t^\op}$, then the functor
\begin{align*}
  \functor{B\tensorcheck_{\Sigma_t}(-)^{\tensorcheck t}}{\SymSeq}{\SymSeq}
\end{align*}
preserves weak equivalences between cofibrant objects, and hence its total left derived functor exists.
\item [(b)] If $Z\in\SymSeq$ is cofibrant, then the functor
\begin{align*}
  \functor{-\tensorcheck_{\Sigma_t} Z^{\tensorcheck t}}{\SymSeq^{\Sigma_t^\op}}{\SymSeq}
\end{align*}
preserves weak equivalences.
\end{itemize}
\end{customprop}

\begin{customprop}{6.11*}
\label{prop:six_point_one_one_star}
Let $t\geq 1$ and consider $\SymSeq$ with the positive flat stable model structure. If $B\in\SymSeq^{\Sigma_t^\op}$, then the functor
\begin{align*}
  \functor{B\tensorcheck_{\Sigma_t}(-)^{\tensorcheck t}}{\SymSeq}{\SymSeq}
\end{align*}
sends cofibrations between cofibrant objects to monomorphisms.
\end{customprop}

All references to Propositions 4.28, 4.29, and 6.11 in the proofs of the main results in \cite{Harper_Spectra} should be replaced by references to Propositions \ref{prop:four_point_two_eight_star}, \ref{prop:four_point_two_nine_star}, and \ref{prop:six_point_one_one_star}, respectively, which are proved below in Section \ref{sec:proofs}.

Propositions 1.6 and 7.7(a) of \cite{Harper_Spectra} are special cases of the statement of Proposition 4.28(a) of \cite{Harper_Spectra}, and hence are incorrect as stated; the following propositions correspond to their corrected versions, respectively, and are special cases of Proposition \ref{prop:four_point_two_eight_star} above.

\begin{customprop}{1.6*}
Let $B\in\bigl(\Spectra\bigr)^{\Sigma_t^\op}$, $t\geq 1$, and $n\geq 0$. If $\function{i}{X}{Y}$ is a cofibration between cofibrant objects in symmetric spectra with the positive flat stable model structure, then the map $B\Smash X^{\wedge t}\rarrow B\Smash Y^{\wedge t}$, after evaluation at $n$, is a cofibration of $\Sigma_t$-diagrams in pointed simplicial sets.
\end{customprop}

\begin{customprop}{7.7*}
Let $B\in\bigl(\Spectra\bigr)^{\Sigma_t^\op}$, $t\geq 1$, and $n\geq 0$. If $\function{i}{X}{Y}$ is a cofibration between cofibrant objects in $\Spectra$ with the positive flat stable model structure, then 
\begin{itemize}
\item [(a)] the map $B\Smash X^{\wedge t}\rarrow B\Smash Y^{\wedge t}$, after evaluation at $n$, is a cofibration in $\sSet_*^{\Sigma_t}$ with the projective model structure inherited from $\sSet_*$,
\item [(b)] the map $B\Smash_{\Sigma_t}Q_{t-1}^t\rarrow B\Smash_{\Sigma_t}Y^{\wedge t}$ is a monomorphism.
\end{itemize}
\end{customprop}

\section{Proofs}
\label{sec:proofs}

The purpose of this section is to prove Propositions \ref{prop:four_point_two_eight_star}, \ref{prop:four_point_two_nine_star}, and \ref{prop:six_point_one_one_star}. The proofs follow closely our original arguments in \cite{Harper_Spectra}.

The following proposition is a useful warm-up for the proof of Proposition \ref{prop:four_point_two_eight_star}.

\begin{prop}
\label{prop:warmup_calculation}
Let $B\in\SymSeq^{\Sigma_t^\op}$, $t\geq 2$, and $r,n\geq 0$. Let $\alpha\geq 1$, $q_0\geq 0$, and $q_1,\dotsc,q_\alpha\geq 1$ such that $q_0+q_1+\dotsb+q_\alpha=t$. If $Z$ is a cofibrant object in $\SymSeq$ with the positive flat stable model structure, then the symmetric sequence
\begin{align*}
  B\tensorcheck
  \Bigl(
  \Sigma_t\cdot_{\Sigma_{q_0}\times\Sigma_{q_1}\times\dotsb\times\Sigma_{q_\alpha}}
  Z^{\tensorcheck q_0}\tensorcheck X_1^{\tensorcheck q_1}\tensorcheck
  \dotsb\tensorcheck X_\alpha^{\tensorcheck q_\alpha}
  \Bigr)
\end{align*}
equipped with the diagonal $\Sigma_t$-action, after evaluation at $[\mathbf{r}]_n$, is a cofibrant object in $\sSet_*^{\Sigma_t}$ with the projective model structure inherited from $\sSet_*$. Here, each $K_i\rarrow L_i$ is a generating cofibration for $\sSet_*$ $(1\leq i\leq\alpha)$, and each $X_i$ is defined as
\begin{align*}
  X_i := G_{p_i}\bigl(S\tensor G_{m_i}^{H_i}(L_i/K_i)\bigr),\quad\quad
  1\leq i\leq \alpha,
\end{align*}
by applying the indicated functors in \cite[4.1]{Harper_Spectra} to the pointed simplicial set $L_i/K_i$, where $m_i\geq 1$, $H_i\subset\Sigma_{m_i}$ a subgroup, and $p_i\geq 0$; in other words, each $X_i$ is assumed to be the cofiber of a generating cofibration for $\SymSeq$ with the positive flat stable model structure.
\end{prop}

\begin{proof}
This is an exercise left to the reader; the argument is by induction on $q_0$, together with (i) the filtrations described in \cite[4.14]{Harper_Spectra} and (ii) the fact that every cofibration of the form $*\rarrow Z$ in $\SymSeq$ is a retract of a (possibly transfinite) composition of pushouts of maps as in \cite[6.17]{Harper_Spectra}, starting with $Z_0=*$.
\end{proof}

\begin{proof}[Proof of Proposition \ref{prop:four_point_two_eight_star}(a)]
Let $m\geq 1$, $H\subsetof \Sigma_m$ a subgroup, and $k,p\geq 0$. Let $\function{g}{\partial\Delta[k]_+}{\Delta[k]_+}$ be a generating cofibration for $\sSet_*$ and consider the pushout diagram \cite[6.17]{Harper_Spectra} in $\SymSeq$ with $Z_0$ cofibrant. It follows from \cite[6.13]{Harper_Spectra} that the diagrams
\begin{align*}
\xymatrix{
  Q_{t-1}^t(g_*)\ar[r]\ar[d] & Q_{t-1}^t(i_0)\ar[d]\\
  D^{\tensorcheck t}\ar[r] & Z_1^{\tensorcheck t}
}\quad\quad
\xymatrix{
  B\tensorcheck Q_{t-1}^t(g_*)\ar[r]\ar[d]^{(*)} & 
  B\tensorcheck Q_{t-1}^t(i_0)\ar[d]^{(**)}\\
  B\tensorcheck D^{\tensorcheck t}\ar[r] & 
  B\tensorcheck Z_1^{\tensorcheck t}
}
\end{align*}
are pushout diagrams in $\SymSeq^{\Sigma_t}$; here, the right-hand diagram is obtained by applying $B\tensorcheck-$ to the left-hand diagram. Since $m\geq 1$, it follows from \cite[3.7]{Harper_Spectra} that $(*)$, after evaluation at $[\mathbf{r}]_n$, is a cofibration in $\sSet_*^{\Sigma_t}$; hence $(**)$, after evaluation at $[\mathbf{r}]_n$, is a cofibration in $\sSet_*^{\Sigma_t}$. Consider a sequence 
\begin{align}
\label{eq:filtration_from_small_obj}
\xymatrix{
  Z_0\ar[r]^{i_0} & Z_1\ar[r]^{i_1} & Z_2\ar[r]^{i_2} & \dotsb\
}
\end{align}
of pushouts of maps as in \cite[6.17]{Harper_Spectra} with $Z_0$ cofibrant, define $Z_\infty:=\colim_q Z_q$, and consider the naturally occurring map $\function{i_\infty}{Z_0}{Z_\infty}$. Using \cite[4.14]{Harper_Spectra} together with Proposition \ref{prop:warmup_calculation}, it is easy to verify that the maps $B\tensorcheck Z_q^{\tensorcheck t}\rarrow B\tensorcheck Q_{t-1}^t(i_q)$ and $B\tensorcheck Q_{t-1}^t(i_q)\rarrow B\tensorcheck Z_{q+1}^{\tensorcheck t}$, after evaluation at $[\mathbf{r}]_n$, are cofibrations in $\sSet_*^{\Sigma_t}$. It follows immediately that each $B\tensorcheck Z_{q}^{\tensorcheck t}\rarrow B\tensorcheck Z_{q+1}^{\tensorcheck t}$, after evaluation at $[\mathbf{r}]_n$, is a cofibration in $\sSet_*^{\Sigma_t}$, and hence the map $B\tensorcheck Z_0^{\tensorcheck t}\rarrow B\tensorcheck Z_\infty^{\tensorcheck t}$, after evaluation at $[\mathbf{r}]_n$, is a cofibration in $\sSet_*^{\Sigma_t}$. Noting that every cofibration between cofibrant objects in $\SymSeq$ with the positive flat stable model structure is a retract of a (possibly transfinite) composition of pushouts of maps as in \cite[6.17]{Harper_Spectra} finishes the proof. 
\end{proof}

The following proposition is an exercise left to the reader.

\begin{prop}
\label{prop:detecting_cofibrations}
Let $G$ be a finite group. Consider any pullback diagram
\begin{align*}
\xymatrix{
  A\ar[d]\ar[r] & C\ar[d]\\
  B\ar[r]^-{f} & D
}
\end{align*}
of monomorphisms in $\sSet_*^G$. If $f$ is a cofibration in $\sSet_*^G$, then the pushout corner map $B\amalg_A C\rarrow D$ is a cofibration in $\sSet_*^G$. 
\end{prop}

\begin{defn}
Let $\I$ be the poset $\{0\rightarrow 1\rightarrow 2\}$, $\I\rarrow\SymSeq$ a diagram, and $t\geq 1$. Consider any subset $\capA\subset\{0\rightarrow 1\rightarrow 2\}^{\times t}=\I^{\times t}$ closed under the canonical $\Sigma_t$-action on $\I^{\times t}$. Denote by
$
  Q_\capA^t:=\colim(\capA\subset\I^{\times t}\rightarrow
  \SymSeq^{\times t}\xrightarrow{\tensorcheck}\SymSeq)
$
the indicated colimit in $\SymSeq$, equipped with the induced $\Sigma_t$-action.
\end{defn}

The following proposition is proved in Pereira \cite{Pereira_spectral_operad}. It provides a refinement of the filtrations for tensor powers of a single map $X\rarrow Y$ in \cite[4.13]{Harper_Spectra} to tensor powers of a composition of maps $X\rarrow Y\rarrow Z$, and will be used in the proof of Proposition \ref{prop:four_point_two_eight_star}(b) below.

\begin{prop}
\label{prop:refined_filtration}
Let $X\xrightarrow{i}Y\xrightarrow{j}Z$ be morphisms in $\SymSeq$ and $t\geq 1$. Consider any convex subset $\capA\subset\{0\rightarrow 1\rightarrow 2\}^{\times t}=\I^{\times t}$ closed under the canonical $\Sigma_t$-action on $\I^{\times t}$. Let $e\in\capA$ be maximal and define
\begin{align*}
  \capA'&:= \capA-\orbit(e),\quad\quad
  \capA_e:= \{v\in\I^{\times t}:\,v\leq e,\quad v\neq e\}.
\end{align*}
Suppose $\capA'\ni(0,\dotsb,0)$. Then $\capA_e\subset\capA'$ and
\begin{itemize}
\item[(a)] the induced map $Q_{\capA'}^t\rarrow Q_\capA^t$ fits into a pushout diagram of the form
\begin{align*}
\xymatrix{
  \Sigma_t\cdot_{\Sigma_p\times\Sigma_q\times\Sigma_r}
  Q_{\capA_e}^t\ar@<6ex>[d]\ar[r] & 
  Q_{\capA'}^t\ar[d]\\
  \quad\quad\quad\quad\quad\
  \Sigma_t\cdot_{\Sigma_p\times\Sigma_q\times\Sigma_r}
  X^{\tensorcheck p}\tensorcheck Y^{\tensorcheck q}\tensorcheck 
  Z^{\tensorcheck r}\ar[r] & 
  Q_\capA^t
} \quad\quad\quad\quad\quad\
\end{align*}
\item[(b)] the induced map $Q_{\capA_e}^t\rarrow X^{\tensorcheck p}\tensorcheck Y^{\tensorcheck q}\tensorcheck Z^{\tensorcheck r}$ is isomorphic to $X^{\tensorcheck p}\tensorcheck -$ applied to the pushout corner map of the commutative diagram
\begin{align*}
\xymatrix{
  Q_{q-1}^q(i)\tensorcheck Q_{r-1}^r(j)\ar[d]_-{\id\tensorcheck j_*}
  \ar[r]^-{i_*\tensorcheck\id} & 
  Y^{\tensorcheck q}\tensorcheck Q_{r-1}^r(j)\ar[d]^-{\id\tensorcheck j_*}\\
  Q_{q-1}^q(i)\tensorcheck Z^{\tensorcheck r}
  \ar[r]_-{i_*\tensorcheck\id} &
  Y^{\tensorcheck q}\tensorcheck Z^{\tensorcheck r}
}
\end{align*}
\end{itemize}
Here, $p:=l_0(e)$, $q:=l_1(e)$, $r:=l_2(e)$, where the ``$i$-length of $e$'', $l_i(e)$, denotes the number of $i$'s in the $t$-tuple $e$, and $Q_{-1}^0:=*$.
\end{prop}

\begin{proof}
This follows from the fact that $\capA_e=\capA_e^1\cup\capA_e^2$ can be written as the union of the convex subsets
\begin{align*}
  \capA_e^1&:=\{v\in\I^{\times t}:\,
  v\leq e,\quad v_j < e_j = 1\quad \text{for some $1\leq j\leq t$}\},\\
  \capA_e^2&:=\{v\in\I^{\times t}:\,
  v\leq e,\quad v_j < e_j = 2\quad \text{for some $1\leq j\leq t$}\}.
\end{align*}
of $\I^{\times t}$, together with the observation in Goodwillie \cite[2.8]{Goodwillie_calc2} that convexity of $\capA_e^1$ and $\capA_e^2$ implies that the commutative diagram
\begin{align*}
\xymatrix{
  \colim_{\capA_e^1\cap\capA_e^2}\capX\ar[r]\ar[d] & 
  \colim_{\capA_e^2}\capX\ar[d]\\
  \colim_{\capA_e^1}\capX\ar[r] & 
  \colim_{\capA_e^1\cup\capA_e^2}\capX
}
\end{align*}
is a pushout diagram in $\SymSeq$, for any functor $\function{\capX}{\I^{\times t}}{\SymSeq}$.
\end{proof}

\begin{rem}
For instance, the induced map $Q_2^3(ji)\rarrow Q_2^3(j)$ is isomorphic to the composition of maps
$
  Q_{\capB_0}^3\rarrow Q^3_{\capB_1}\rarrow Q^3_{\capB_2}
  \rarrow Q^3_{\capB_3}
$
where
\begin{align*}
  \capB_0&:=\{v\in\I^{\times 3}:\,l_0(v)\geq 1\},&\quad&\capB_1:=\capB_0\cup\orbit\bigl((1,1,1)\bigr),\\
  \capB_2&:=\capB_1\cup\orbit\bigl((1,1,2)\bigr),&\quad&
  \capB_3:=\capB_2\cup\orbit\bigl((1,2,2)\bigr).
\end{align*}
\end{rem}

\begin{proof}[Proof of Proposition \ref{prop:four_point_two_eight_star}(b)]
Proceed as above for part (a) and consider the commutative diagram 
\begin{align}
\label{eq:ladder_diagram}
\xymatrix{
  B\tensorcheck Z_0^{\tensorcheck t}\ar@{=}[d]\ar[r] & 
  B\tensorcheck Q_{t-1}^t(i_0)\ar[d]\ar[r] &
  B\tensorcheck Q_{t-1}^t(i_1 i_0)\ar[d]\ar[r] &
  B\tensorcheck Q_{t-1}^t(i_2 i_1 i_0)\ar[d]\ar[r] &
  \dotsb\\
  B\tensorcheck Z_0^{\tensorcheck t}\ar[r] &
  B\tensorcheck Z_1^{\tensorcheck t}\ar[r] &
  B\tensorcheck Z_2^{\tensorcheck t}\ar[r] &
  B\tensorcheck Z_3^{\tensorcheck t}\ar[r] &
  \dotsb
}
\end{align}
in $\SymSeq^{\Sigma_t}$. We know by part (a) that the bottom row, after evaluation at $[\mathbf{r}]_n$, is a diagram of cofibrations in $\sSet_*^{\Sigma_t}$. Using Propositions \ref{prop:refined_filtration}, \ref{prop:detecting_cofibrations}, and \ref{prop:warmup_calculation}, together with \cite[4.14]{Harper_Spectra}, it is easy to verify that each of the maps
\begin{align*}
  &B\tensorcheck Q_{t-1}^t(i_0)\rarrow 
  B\tensorcheck Z_1^{\tensorcheck t},\\
  &B\tensorcheck Q_{t-1}^t(i_1 i_0)\rarrow 
  B\tensorcheck Q_{t-1}^t(i_1)\rarrow
  B\tensorcheck Z_2^{\tensorcheck t},\\
  &B\tensorcheck Q_{t-1}^t(i_2 i_1 i_0)\rarrow 
  B\tensorcheck Q_{t-1}^t(i_2 i_1)\rarrow
  B\tensorcheck Q_{t-1}^t(i_2)\rarrow
  B\tensorcheck Z_3^{\tensorcheck t},\quad\dotsb
\end{align*}
and hence the vertical maps in \eqref{eq:ladder_diagram},
after evaluation at $[\mathbf{r}]_n$, are cofibrations in $\sSet_*^{\Sigma_t}$. It follows that applying $\colim_{\Sigma_t}(-)$ to \eqref{eq:ladder_diagram} gives the commutative diagram \cite[6.20]{Harper_Spectra} of monomorphisms, hence the induced map $B\tensorcheck_{\Sigma_t}Q_{t-1}^t(i_\infty)\rarrow B\tensorcheck_{\Sigma_t}Z_\infty^{\tensorcheck t}$ is a monomorphism. Noting that every cofibration between cofibrant objects in $\SymSeq$ is a retract of a (possibly transfinite) composition of pushouts of maps as in \cite[6.17]{Harper_Spectra}, together with \cite[6.14]{Harper_Spectra}, finishes the proof.
\end{proof}

The following proposition, which appeared in an early version of \cite{Schwede_book_project}, can be thought of as a refinement of the arguments in \cite[15.5]{Mandell_May_Schwede_Shipley} and \cite[3.3]{Shipley_comm_ring}. 
\begin{prop}
\label{prop:connectivity_of_G_orbits}
Let $G$ be a finite group, $Z'\rarrow Z$ a morphism in $(\Spectra)^G$, and $k\in\Z\cup\{\infty\}$. Assume that $G$ acts freely on $Z',Z$ away from the basepoint $*$, and consider the $G$-orbits spectrum $Z/G:=\colim_G Z\Iso S\Smash_G Z$. If $Z$ (resp. $Z'\rarrow Z$) is $k$-connected, then $Z/G$ (resp. $Z'/G\rarrow Z/G$) is $k$-connected.
\end{prop}

\begin{proof}
Consider the contractible simplicial set $EG\xrightarrow{\wequiv}*$ with free right $G$-action, given by realization of the usual simplicial bar construction with respect to Cartesian product $EG=|\BAR^\times(*,G,G)|$. Since $G$ acts freely on $Z$ away from the basepoint,  the induced map $EG_+\Smash_G Z\xrightarrow{\wequiv}*_+\Smash_G Z\Iso S\Smash_G Z$ of symmetric spectra is a weak equivalence. We need to verify that $S\Smash_G Z$ is $k$-connected; it suffices to verify that $EG_+\Smash_G Z$ is $k$-connected. The symmetric spectrum $EG_+\Smash_G Z$ is isomorphic to the realization of the usual simplicial bar construction with respect to smash product $|\BAR^\wedge(*_+,G_+,Z)|$. We know by assumption that $Z$ is $k$-connected, hence $\BAR^\wedge(*_+,G_+,Z)$ is objectwise $k$-connected. The other case is similar.
\end{proof}

\begin{proof}[Proof of Proposition \ref{prop:four_point_two_nine_star}]
Consider part (b). Suppose $A\rarrow B$ in $\SymSeq^{\Sigma_t^\op}$ is a weak equivalence. Then it follows from Propositions \ref{prop:four_point_two_eight_star}(a) and \ref{prop:connectivity_of_G_orbits} (with $k=\infty$) that the induced map $A\tensorcheck_{\Sigma_t} Z^{\tensorcheck t}\rarrow B\tensorcheck_{\Sigma_t} Z^{\tensorcheck t}$ is a weak equivalence. Consider part (a). Suppose $X\rarrow Y$ in $\SymSeq$ is a weak equivalence between cofibrant objects; we want to show that $B\tensorcheck_{\Sigma_t} X^{\tensorcheck t}\rarrow B\tensorcheck_{\Sigma_t} Y^{\tensorcheck t}$ is a weak equivalence. The map $*\rarrow B$ factors in $\SymSeq^{\Sigma_t^\op}$ as 
$*\rarrow B^c\rarrow B$ a cofibration followed by an acyclic fibration, 
\begin{align}
\label{eq:commutative_square_fatten_up_B}
\xymatrix{
  B^c\tensorcheck_{\Sigma_t} X^{\tensorcheck t}\ar[r]\ar[d] & 
  B^c\tensorcheck_{\Sigma_t} Y^{\tensorcheck t}\ar[d]\\
  B\tensorcheck_{\Sigma_t} X^{\tensorcheck t}\ar[r] & 
  B\tensorcheck_{\Sigma_t} Y^{\tensorcheck t}
}
\end{align}
diagram \eqref{eq:commutative_square_fatten_up_B} commutes, and since three of the maps are weak equivalences, so is the fourth; here, we have used \cite[4.29(b)]{Harper_Spectra}.
\end{proof}

\begin{proof}[Proof of Proposition \ref{prop:six_point_one_one_star}]
Suppose $X\rarrow Y$ in $\SymSeq$ is a cofibration between cofibrant objects; we want to show that $B\tensorcheck_{\Sigma_t} X^{\tensorcheck t}\rarrow B\tensorcheck_{\Sigma_t} Y^{\tensorcheck t}$ is a monomorphism. This follows immediately from Proposition \ref{prop:four_point_two_eight_star}.
\end{proof}

\bibliographystyle{plain}
\bibliography{CorrigendumHomotopyModulesSpectra.bib}

\end{document}